\documentclass[11pt,psamsfonts,oneside]{amsart}
\usepackage{latexsym}
\usepackage{amsmath, amsfonts, amssymb, amscd}
\usepackage[all]{xy}
\usepackage[height=22.5cm, width=16cm]{geometry}
\usepackage[latin1]{inputenc}
\usepackage[colorlinks=true, pdfstartview=FitV, linkcolor=blue, citecolor=blue, urlcolor=blue]{hyperref}

\newtheorem{definition}{Definition}[section]
\newtheorem{theorem}[definition]{Theorem}
\newtheorem{lemma}[definition]{Lemma}
\newtheorem{corollary}[definition]{Corollary}
\newtheorem{proposition}[definition]{Proposition}
\newtheorem*{mainTheorem}{Theorem 1}

\theoremstyle{definition}

\theoremstyle{remark}
\newtheorem*{remark}{Remark}

\DeclareMathOperator{\Ad}{Ad} 
 
\DeclareMathOperator{\codim}{codim}

\newcommand{\im}{\mathrm{i}}
\newcommand{\ce}{\mathbb{C}}
\newcommand{\er}{\mathbb{R}}

\newcommand{\quot}[2]{#1/\!\!/#2}
\newcommand{\Mp}{\mathcal M_\liep}
\newcommand{\Mq}{\mathcal M_\lieq}
\newcommand{\sgmp}{\mathcal S_G(\Mp)}
\newcommand{\sgmpb}{\mathcal S_G(\Mp(\beta))}

\newcommand{\lieg}{\mathfrak{g}}
\newcommand{\liez}{\mathfrak{z}}

\newcommand{\liep}{\mathfrak{p}}
\newcommand{\liek}{\mathfrak{k}}
\newcommand{\lieu}{\mathfrak{u}}

\newcommand{\lieq}{\mathfrak{q}}

\newcommand{\pw}{\mathbb P(W)}
\newcommand{\pv}{\mathbb P(V)}

\newcommand{\phatv}{\mathbb P(\hat V)}

\newcommand{\ddt}{\left.\frac{\partial}{\partial
t}\right|_0}

\title{Closed orbits of real reductive representations}
\author{Henrik St\"otzel}

\address{Fakult\"at f\"ur Mathematik\\
  Ruhr Universit\"at Bochum\\
  Universit\"atsstrasse 150\\
  D - 44780 Bochum}
\email{henrik.stoetzel@rub.de}
\thanks{The author is supported by the
  Sonderforschungsbereich SFB/TR12 of the Deutsche
  Forschungsgemeinschaft}

\begin{document}

\begin{abstract}
We prove that the set of closed orbits in a real reductive representation contains a subset which is open with respect to the real Zariski topology if it has non-empty interior. In particular the set of closed orbits is dense.
\end{abstract}

\maketitle

\section{Introduction}

We consider the representation $G\times V\to V$ of a real Lie group on a finite dimensional real vector space $V$.
We say that the representation is real reductive, if there exists a complex reductive Lie group $U^\ce$ with Cartan involution $\theta$ such that $G$ is a $\theta$-stable subgroup of $U^\ce$ with finitely many connected components and if the representation $G\times V\to V$ extends to a holomorphic representation of $U^\ce$ on a complex vector space containing $V$ as a $G$-invariant subspace. 

Denoting by $V_c$ the set of closed $G$-orbits in $V$,
our main result is the following.

\begin{mainTheorem}\label{theorem:closedOrbits}
Let $G\times V\to V$ be a real reductive representation and $W$ a subspace of $V$ such that $V_c\cap W$ has non-empty interior in $W$. Then $V_c\cap W$ contains a subset $\mathcal U$ of $W$ which is open in $W$ with respect to the real Zariski topology. In particular, $V_c\cap W$ is dense in $W$. If $W$ is $G$-invariant, then $\mathcal U$ can be chosen to be $G$-invariant.
\end{mainTheorem}

Theorem~\ref{theorem:closedOrbits} is a generalization of a result of Birkes for representations of real forms. If $G$ is a real form of $U^\ce$ and if the $G$-representation $V$ is given by restriction of a holomorphic representation of $U^\ce$ on $V^\ce$, then it is shown in \cite{Birkes} that the closed $G$-orbits in $V$ correspond to closed $U^\ce$-orbits in $V^\ce$. We will give the usual proof using moment maps for this result. Using a stratification of $V$ defined by orbit types of closed orbits, this enables us to prove Theorem~\ref{theorem:closedOrbits} for real forms. 

In our general situation, where $G$ is not a real form, the idea is to complexify the $G$-representation $V$. If $G$ is semisimple, $G$ is a compatible real form of its universal complexification $G^\ce$ and the $G$-representation $V$ extends to a holomorphic $G^\ce$-representation $V^\ce$, so the proof for real forms applies. If $G$ is commutative, the weight space decomposition can be used in order to give an explicit description of the set of closed orbits.

The mixed case, where $G$ is a product of its center and a semisimple subgroup, is more sophisticated. We can not describe the closed $G$-orbits in terms of the closed orbits of the semisimple subgroup and the center. Instead, we will construct a group $H$ which contains $G$, which differs from $G$ only by a subgroup of its center and which is a real form of a complex reductive group $U^\ce$. Moreover $H$ admits a representation on $V$ which extends the $G$-representation and which is given by restriction of a holomorphic $U^\ce$-representation on $V^\ce$. Thus Theorem~\ref{theorem:closedOrbits} for real forms applies to the representation of $H$. Then it remains to establish a relation between the closed $G$-orbits and the closed $H$-orbits. Our main tool is the fact, that the closed $G$-orbits inside the set of semistable points $\mathcal S_G(\Mp):=\{z\in Z; \overline{G\cdot z}\cap\mu_\liep^{-1}(0)\neq\emptyset\}$ with respect to a restricted moment map $\mu_\liep$ (see Section~\ref{section:rep}) are parametrized by the orbits of a maximal compact subgroup of $G$ in the zero-fiber $\mu_\liep^{-1}(0)$ (\cite{HSch}, \cite{Kirwan}, \cite{Ness}). In our situation a $G$-orbit is closed inside the set of semistable points with respect to a standard moment map if the corresponding $H$-orbit is closed inside the set of semistable points with respect to a shifted standard moment map. Since we do not know if the latter set of semistable points is dense, we give an equivalent formulation of Theorem~\ref{theorem:closedOrbits} in projective space. Here the set of semistable points is dense by \cite{HMig}.

It should be noted that in the special case where $G$ is complex reductive and where its representation on $V$ is holomorphic, the situation is much simpler. Here the strata defined with respect to orbit types of closed orbits are locally closed with respect to the complex Zariski topology. Then the assumption that the set of closed orbits has non-empty interior in $W$ implies that $\mathcal U$ can be chosen to be the intersection of a stratum with $W$. In particular all orbits which intersect $\mathcal U$ are of the same orbit type. In the real case this is wrong. For this consider the adjoint representation of a real semisimple group. Here there exist finitely many Cartan subalgebras. An orbit through a regular element is closed and the isotropy equals the centralizer of the Cartan subalgebra. Therefore each Cartan subalgebra defines an open stratum which consists of closed orbits and the union of these open strata is dense.

Theorem~\ref{theorem:closedOrbits} may also be applied to actions on Kählerian manifolds $Z$ which admit a restricted moment map $\mu_\liep$ with respect to the action of $G$. If $z\in Z$ is contained in the zero fiber $\mu_\liep^{-1}(0)$, then there exists a geometric $G$-slice at $z$ (\cite{HSch}). In particular the $G$-action in a neighborhood of $z$ is described by the isotropy representation of $G_z$ on the tangent space $T_zZ$. Applying Theorem~\ref{theorem:closedOrbits} to the isotropy representation, we get

\begin{corollary}\label{corollary:manifolds}
Assume that there exists a restricted moment map $\mu_\liep$ for the action of $G$ on $Z$ and that the topological Hilbert quotient $\quot {\mathcal S_G(\Mp)}G$ is connected. If the set of closed $G$-orbits inside the set of semistable points $\mathcal S_G(\Mp)$ has non-empty interior, then its interior is dense in $\mathcal S_G(\Mp)$.
\end{corollary}

An application of Theorem~\ref{theorem:closedOrbits} is given in \cite{St}. If $H$ is a generic isotropy of the $G$-action on $Z$, define $Y:=\{z\in Z; G_z=H\text{ and $G\cdot z$ is closed}\}$ and let $\overline Y$ denote the closure of $Y$. Let $\mathcal N_G(H)$ denote the normalizer of $H$ in $G$. The main result in \cite{St} is a description of the topological Hilbert quotient $\quot ZG$ in terms of the quotient $\quot Y{\mathcal N_G(H)}$. Here Theorem~\ref{theorem:closedOrbits} applies to show that $Y$ is smooth.

The author would like to thank P.~Heinzner for
carefully reading the manuscript of this paper.

\section{Real reductive representations}

Let $U$ be a compact Lie group and $U^\ce$ its complexification. We have a Cartan decomposition $U^\ce=U\exp(\im\lieu)$ where $\lieu$ is the Lie algebra of $U$. The corresponding Cartan involution $\theta$ is given by $\theta(u\exp(\xi))=u\exp(-\xi)$. We call a closed real subgroup $G$ of $U^\ce$ \emph{compatible}, if it is invariant under the Cartan involution $\theta$ and consists of only finitely many connected components. Equivalently, a Lie subgroup $G$ is compatible if and only if $G=K\exp(\liep)$ for a compact subgroup $K$ of $U$ and a subspace $\liep$ of $\im\lieu$. The decomposition $G=K\exp(\liep)$ is called the \emph{Cartan decomposition} of $G$.
Note that in particular $U^\ce$ is a compatible subgroup of $U^\ce$ with Cartan decomposition $U^\ce=U\exp(\im\lieu)$ where $\lieu$ denotes the Lie algebra of $U$.

We call a representation $G\times V\to V$ of a real Lie group on a real vector space a \emph{real reductive representation}, if $G$ is a compatible subgroup of a complex reductive group $U^\ce$ and if $V$ is a $G$-invariant real subspace of a holomorphic $U^\ce$-representation space $\hat V$.

\section{Semistability in projective space}\label{section:rep}

The linear action $U^\ce\times\hat V\to\hat V$ induces a holomorphic action of $U^\ce$ on the complex projective space $\phatv$ such that the canonical projection $\pi\colon \hat V\setminus\{0\}\to\phatv$ is $U^\ce$-equivariant, i.\,e.~commutes with the actions of $U^\ce$ on $\hat V$ and on $\phatv$.
Our goal is to give an equivalent formulation of Theorem~\ref{theorem:closedOrbits} in terms of the induced $G$-action on projective space. For this, we describe the set of closed $G$-orbits in $V$ in terms of a moment map on $\phatv$.

We fix a Hermitian inner product on $\hat V$ such that $U$ acts by unitary operators on $\hat V$. Recall that the associated Fubini Study metric $\omega$ is a Kählerian metric on $\phatv$, which is invariant under the action of the unitary group and which is in particular invariant under the action of $U$. The corresponding moment map $\mu\colon \phatv\to\lieu^*$ is defined by the equation $d\mu^\xi=\iota_{\xi_*}\omega$. Here $\xi_*(x)=\ddt\exp(t\xi)\cdot x$ is the vectorfield induced by the action of the one-parameter group $\{\exp(t\xi);t\in\er\}$ and $\iota_{\xi_*}$ is the contraction of $\omega$ with $\xi_*$. Moreover, we require $\mu$ to be equivariant with respect to the coadjoint action of $U$ on $\lieu^*$ which is given by $(g\cdot \varphi)(\xi)=\varphi(\Ad(g)\xi)$.
Explicitly a moment map is given by $\mu^\xi(\pi(v))=\im\frac{\langle \xi_*v,v\rangle}{\langle v,v\rangle}$. It is unique up to addition of an element in the dual of the center of $\lieu^*$. The moment map contains much information on the geometry of the $U^\ce$-action on $\phatv$ (\cite{Kirwan}, \cite{Ness}). Analogously the geometry of the $G$-action, which we are interested in, can be described in terms of the restriction of the moment map to $(\im\liep)^*$ (\cite{HSch}). In order to simplify notation, we identify $(\im\liep)^*$ with $\liep$. Then the \emph{restricted moment map} $\mu_\liep\colon\phatv\to\liep$ is given by $\mu_\liep^\xi(x)=\langle\mu_\liep(x),\xi\rangle=\mu^{-\im\xi}(x)=\frac{\langle \xi_*v,v\rangle}{\langle v,v\rangle}$

For $\beta\in\liep$, we define $\Mp(\beta):=\mu_\liep^{-1}(\beta)$ and in order to shorten notation we set $\Mp:=\mu_\liep^{-1}(0)$. Moreover we define $\mathcal S_G(\Mp(\beta)):=\{x\in\phatv;\overline{G\cdot x}\cap\Mp(\beta)\neq\emptyset\}$ to be the set of $G$-orbits which intersect $\Mp(\beta)$ in their closure. For $\beta=0$ we call $\mathcal S_G(\Mp):=\mathcal S_G(\Mp(0))$ the set of \emph{semistable points} in $\phatv$ with respect to $\mu_\liep$. If $\beta\in\liep$ is contained in the center of $\lieu^\ce$, the shifted moment map $\mu+\im\beta$ is again a moment map and $\Mp(\beta)$ is the zero-fiber of the restricted moment map $\mu_\liep-\beta$, so $\sgmpb$ is the set of semistable points with respect to $\mu_\liep-\beta$. We have a relation between the closed orbits in $\mathcal S_G(\Mp(\beta))$ and $\Mp(\beta)$. 

\begin{theorem}[\cite{HSch}]\label{theorem:semistable}
Let $\beta\in\liep$ be contained in the center of $\lieu^\ce$.
\begin{enumerate}\item A $G$-orbit $G\cdot x\subset\sgmpb$ is closed in $\sgmpb$ if and only if $G\cdot x\cap\Mp(\beta)\neq\emptyset$.\item Every non-closed orbit in $\sgmpb$ contains a unique closed orbit in its closure and has strictly larger dimension than that closed orbit.
\end{enumerate}
\end{theorem}

In particular, the closed $G$-orbits in the set of semistable points $\sgmp$ are exactly those orbits which intersect the zero-fiber $\Mp=\{\pi(v);<\xi_*v,v>=0\text{ for all }\xi\in\liep\}$ of $\mu_\liep$. It is shown in \cite{RS} that a $G$-orbit $G\cdot v$ in $\hat V$ is closed if and only if it intersects $\mathcal M:=\{v\in\hat V;<\xi_*v,v>=0\text{ for all }\xi\in\liep\}$. In our terminology, $\mathcal M$ is the zero-fiber of a restricted moment map on $\hat V$ for which the set of semistable points equals $\hat V$. We observe that $\mathcal M\setminus\{0\}=\pi^{-1}(\Mp)$. Since $\pi$ is $G$-equivariant, together with Theorem~\ref{theorem:semistable} this proves

\begin{lemma}\label{lemma:ClosedObenUnten}
An orbit $G\cdot v$, $v\neq 0$ is closed in $\hat V$ if and only if $G\cdot \pi(v)$ is closed in $\sgmp$.
\end{lemma}

The set of closed $G$-orbits in $\hat V$ is $\ce^*$-invariant. Therefore, applying Lemma~\ref{lemma:ClosedObenUnten}, Theorem~\ref{theorem:closedOrbits} can be reformulated as follows.

\medskip

\emph{Let $W$ be a linear subspace of $V$ and assume that the intersection of the set of closed orbits in $\sgmp$ with $\pi(W)$ has non-empty interior in $\pi(W)$. Then it contains a subset $\mathcal U$ of $\pi(W)$ which is open in $\pi(W)$ with respect to the real Zariski topology. If $\pi(W)$ is $G$-invariant, then $\mathcal U$ can be chosen to be $G$-invariant.}

\medskip

Now we give an explicit description of the set of semistable points. It is shown in \cite{RS} that each $G$-orbit in $\hat V$ contains exactly one closed $G$-orbit in its closure (compare Theorem~\ref{theorem:semistable}). We define $\mathcal N:=\{v\in \hat V;0\in\overline{G\cdot v}\}$ to be the nullcone of the $G$-representation $\hat V$, i.\,e.~the set of $G$-orbits which contain $0$ as the unique closed orbit in their closure. Note that $\mathcal N$ is $\ce^*$-stable. By a result in \cite{HSchuetz}, the nullcone $\mathcal N$ is a real algebraic subset of $\hat V$. Then $\pi(\mathcal N\setminus\{0\})$ is real algebraic in $\phatv$.

\begin{lemma}\label{lemma:semistabilZOffen}
The set of semistable points is given by $\sgmp=\pi(\hat V\setminus\mathcal N)$. In particular, it is open in $\phatv$ with respect to the real Zariski topology.
\end{lemma}

\begin{proof}
Let $x\in\sgmp\cap\pi(\mathcal N\setminus\{0\})$ and $y\in\overline{G\cdot x}\cap\Mp$. Then $y\in\pi(\mathcal N)$ since $\pi(\mathcal N)$ is closed and $G$-invariant. But $\Mp=\pi(\mathcal M\setminus\{0\})$ and $\mathcal M\cap\mathcal N=\{0\}$, a contradiction. This shows $\sgmp\subset\pi(\hat V\setminus\mathcal N)$.

Conversely, for $v\in\hat V\setminus\mathcal N$, there exists a $w\in\overline{G\cdot v}\cap(\mathcal M\setminus\{0\})$. Then $\pi(w)\in\overline{G\cdot\pi(v)}\cap\Mp$ by continuity of $\pi$ and thus $\pi(v)\in\sgmp$.
\end{proof}

There is no analogous statement for $\sgmpb$. But if $G$ is complex reductive, the following is known (\cite{HMig}).

\begin{proposition}\label{proposition:SemistabilBetaZOffen}
Assume $G=U^\ce$. Then $\sgmpb=\mathcal S_{U^\ce}(\mathcal M_{\im\lieu}(\beta))$ is open in $\phatv$ with respect to the complex Zariski topology.
\end{proposition}

\section{Actions of real forms on projective space}\label{section:realForms}

We call $G$ a real form of $U^\ce$ if its Lie algebra is a real form of $\lieu^\ce$ and if $G$ intersects every connected component of $U^\ce$. If $G$ is a real form of $U^\ce$ and if $V$ is a real form of $\hat V$, the sets of semistable points and the sets of closed orbits with respects to the actions of $G$ and $U^\ce$ are related. Observe that $\pi(V\setminus\{0\})\subset\phatv$ can be identified with the real projective space $\pv$ since $V$ is totally real in $\hat V$, i.\,e.~$V\cap\im V=\{0\}$ holds.

In order to shorten notation, we define $\mathcal S:=\mathcal S_{U^\ce}(\mathcal M_{\im\lieu}(\beta))\subset\phatv$. We denote by $\sgmpb_c$ and $\mathcal S_c$ the set of closed $G$-orbits and the set of closed $U^\ce$-orbits in $\sgmpb$ and $\mathcal S$, respectively. The following proposition bases on a more general argument which states that the moment map $\mu$ and the restricted moment map $\mu_\liep$ coincide up to a constant on $K$-stable Lagrangian submanifolds.

\begin{proposition}\label{proposition:semistabilGUce}
Let $\beta\in\liep$ be contained in the center of $\lieu^\ce$ and assume that $G$ is a real form of $U^\ce$ and that $V$ is a real form of $\hat V$.
\begin{enumerate}
\item
$\mathcal S_G(\Mp(\beta))\cap\pv=\mathcal S\cap\pv$.
\item $\sgmpb_c\cap\pv=\mathcal S_c\cap\pv$, i.\,e.~a $G$-orbit $G\cdot x$ in $\mathcal S_G(\Mp(\beta))\cap\pv$ is closed if and only if $U^\ce\cdot x$ is closed in $\mathcal S$.
\end{enumerate}
\end{proposition}

\begin{proof}
Recall that we are given a $U$-invariant Hermitian inner product $\langle \cdot,\cdot\rangle$ on $\hat V$. By \cite{RS} this inner product can be chosen such that $V$ is Lagrangian with respect to the Kähler structure which is given by its imaginary part. Since $K$ operates on $V$, this gives $\langle \xi_* v,v\rangle=0$ for $v\in V$ and $\xi\in\liek$. Therefore $\mu^\xi\equiv 0$ for $\xi\in\liek$. Since $G$ is a real form of $U^\ce$, we have a Lie algebra decomposition $\lieu=\liek\oplus\im\liep$ and we conclude $\Mp\cap\pv=\mathcal M_{\im\lieu}\cap\pv$. Since $\beta$ is contained in $\liep$ this also gives $\Mp(\beta)\cap\pv=\mathcal M_{\im\lieu}(\beta)\cap\pv$. 

For the first part of the proposition, the inclusion $\mathcal S\subset \mathcal S_G(\Mp(\beta))$ is shown in \cite{HSch}. Conversely, for $x\in\mathcal S_G(\Mp(\beta))\cap\pv$, the closure $\overline{G\cdot x}$ intersects $\Mp(\beta)\cap\pv=\mathcal M_{\im\lieu}(\beta)\cap\pv$. In particular, $\overline{U^\ce\cdot x}$ intersects $\mathcal M_{\im\lieu}(\beta)\cap\pv$ and thus $x\in \mathcal S$.

For the second part, if $G\cdot x$ is closed in $\mathcal S_G(\Mp(\beta))$, then it intersects $\Mp(\beta)$ by Theorem~\ref{theorem:semistable}. But since $\Mp(\beta)\cap\pv=\mathcal M_{\im\lieu}(\beta)\cap\pv$, it follows that $U^\ce\cdot x$ is closed, again with Theorem~\ref{theorem:semistable}. Conversely, assume $U^\ce\cdot x$ is closed. The action of $U^\ce$ is holomorphic, $G$ is a real form of $U^\ce$ and $\pv$ is totally real in $\phatv$. Therefore the real dimension of each $G$-orbit in $U^\ce\cdot x\cap\pv$ equals the complex dimension of $U^\ce\cdot x$. In particular, all the $G$-orbits in the intersection have the same dimension. Moreover they are contained in $\sgmpb$ by the first part of the proposition. Then it follows from Theorem~\ref{theorem:semistable} that all these $G$-orbits, and in particular $G\cdot x$, are closed.
\end{proof}

\begin{remark}
The same proof applies in order to show the result of Birkes (\cite{Birkes}) that a $G$-orbit $G\cdot v$ in $V$ is closed if and only if $U^\ce\cdot v$ is closed.
\end{remark}

Let $\beta\in\liep$ be contained in the center of $\lieu^\ce$ and consider the shifted moment map $\mu+\im\beta$.
For its set of semistable points $\mathcal S\subset\phatv$, the GIT-Quotient $\pi\colon\mathcal S\to\quot{\mathcal S}{U^\ce}$ in the sense of Mumford (\cite{MFK}) exists and is an affine map. Each fiber of $\pi$ contains a unique closed orbit and every other orbit in the fiber contains the closed orbit in its closure. The quotient defines an orbit type stratification of $\mathcal S$ as follows (\cite{Lu73},\cite{Ri72},\cite{LR79}). For a subgroup $H$ of $U^\ce$, we denote by $\mathcal S^{<H>}:=\{x\in\mathcal S; U^\ce\cdot x\text{ closed, }(U^\ce)_x=H\}$ the set of points on closed orbits, for which the isotropy group is given by $H$. Then the saturation $I_H:=\pi^{-1}(\pi(\mathcal S^{<H>})):=\{x\in\mathcal S; \overline{U^\ce\cdot x}\cap \mathcal S^{<H>}\neq\emptyset\}$ is by definition the \emph{$H$-isotropy stratum}. The closed orbits in $I_H$ are exactly those which intersect $\mathcal S^{<H>}$ and they are the orbits of minimal dimension in $I_H$. The isotropy groups of two points on an orbit are conjugate. Therefore $I_H=I_{gHg^{-1}}$ for $g\in U^\ce$. Since each orbit in $\mathcal S$ contains a closed orbit in its closure, the isotropy strata cover $\mathcal S$. Moreover, two strata $I_{H}$ and $I_{H'}$ are either disjoint or equal, where the latter is the case if and only if $H$ and $H'$ are conjugate. Considering conjugacy classes of isotropy groups,  $\mathcal S$ is the disjoint union of the strata and the union is locally finite. A stratum $I_H$ is open in its closure with respect to the complex Zariski topology and if the closure $\overline{I_H}$ intersects a stratum $I_{H'}$ with $I_H\neq I_{H'}$, then after conjugation, $H$ is contained in $H'$, i.\,e.~there exists a $g\in U^\ce$ such that $gHg^{-1}< H'$. This implies that for $n\in\mathbb N$ the set $\mathcal I_n:=\bigcup_{\codim H\geq n} I_H$ is open in $\mathcal S$ with respect to the complex Zariski topology.

For $n\in\mathbb N$, we define $\mathcal O_n:=\{x\in\phatv; \dim U^\ce\cdot x\geq n\}$ to be the set of orbits with complex dimension at least $n$. Since the elements of the Lie algebra of $U^\ce$ define holomorphic vectorfields on $\phatv$, the set $\mathcal O_n$ is open with respect to the complex Zariski topology. 

Now, for actions of real forms, we prove a result which is slightly more general than the statement of Theorem~\ref{theorem:closedOrbits} since we do not assume $\beta=0$. We need this generality in order to prove Theorem~\ref{theorem:closedOrbits}. Recall that $\mathcal S_G(\mathcal M_\liep(\beta))_c$ denotes the set of closed $G$-orbits in $\mathcal S_G(\mathcal M_\liep(\beta))$.

\begin{proposition}\label{proposition:Projektiv}
Assume $G$ is a real form of $U^\ce$ and $\hat V=V^\ce$. Let $W$ be a subspace of $V$ and let $n$ be maximal such that $\mathcal O_n\cap\pw$ is non-empty. Then, for a $\beta\in\liep$ which is contained in the center of $\lieu^\ce$, the intersection $\mathcal S_G(\mathcal M_\liep(\beta))_c\cap\mathcal O_n\cap\pw$ is open in $\pw$ with respect to the real Zariski topology. 
\end{proposition}

\begin{proof}
For $x\in\pw\subset\pv$, the orbit $G\cdot x$ is closed in $\mathcal S_G(\mathcal M_\liep(\beta))$ if and only if $U^\ce\cdot x$ is closed in $\mathcal S$ (Proposition~\ref{proposition:SemistabilBetaZOffen}). By definition of $n$ and $\mathcal O_n$, an orbit $U^\ce\cdot x$ with $x\in\mathcal O_n\cap\pw$ is of complex dimension $n$. So for $x\in\mathcal S_G(\mathcal M_\liep(\beta))_c\cap\mathcal O_n\cap\pw$ the orbit $U^\ce\cdot x$ is closed and of dimension $n$. Then we have $x\in I_H$ where $H$ is the isotropy group at $x$ and since the dimension of $U^\ce\cdot x$ equals the codimension of the isotropy group, $U^\ce\cdot x$ is contained in the union of isotropy strata $\mathcal I_n=\bigcup_{\codim H\geq n} I_H$.

Conversely, for $x\in \mathcal I_n\cap\mathcal O_n\cap\pw$ the orbit $U^\ce\cdot x$ is of dimension $n$ by the choice of $n$, so it is an orbit of minimal dimension in $\mathcal I_n$. Since the closure of $U^\ce\cdot x$ is contained in $\mathcal I_n$, it follows from Theorem~\ref{theorem:semistable} that $U^\ce\cdot x$ is closed. This shows \[\mathcal S_G(\mathcal M_\liep(\beta))_c\cap\mathcal O_n\cap\pw=\mathcal I_n\cap\mathcal O_n\cap\pw.\]
But $\mathcal I_n$ and $\mathcal O_n$ are open in $\phatv$ with respect to the complex Zariski topology, so the intersection $\mathcal I_n\cap\mathcal O_n\cap\pw$ is open in $\pw$ with respect to the real Zariski topology.
\end{proof}

\section{Proof of the main result}\label{section:enlarging}

We observed in Section~\ref{section:rep} that for the proof of Theorem~\ref{theorem:closedOrbits} it suffices to show that the intersection $\sgmp_c\cap\pi(W)$ of the set of closed orbits inside the set of semistable points with $\pi(W)$ contains a subset which is open in $\pi(W)$ with respect to the real Zariski topology. In Section~\ref{section:realForms} we proved a slightly more general result for actions of real forms. For the general case we will now give a construction of a representation of a real form which contains $G$ in a convenient way. This enables us to prove Theorem~\ref{theorem:closedOrbits} using Proposition~\ref{proposition:Projektiv}.

For our construction, we will make several assumptions on $G$ which do not affect the set $V_c$ of closed orbits.

First, replacing $G$ by its connected component of the identity $G^\circ$, we may assume that $G$ is connected. For this note that $G^\circ$ is compatible with Cartan decomposition $G^\circ=K^\circ\exp(\liep)$ and that a $G$-orbit $G\cdot v$ is closed if and only if $G^\circ\cdot v$ is closed.

The connected compatible subgroup $G$ of $U^\ce$ admits a decomposition $G=G_S\cdot G_Z$, where $G_S$ is a semisimple subgroup of $G$ and $G_Z$ is the center. Here the Lie algebra of $G_S$ is given by the Lie bracket $[\lieg,\lieg]$. Both, $G_S$ and $G_Z$ are compatible subgroups of $U^\ce$. Let $G_Z=K_Z\exp(\liep_Z)$ denote the Cartan decomposition of the center and define $P_Z:=\exp(\liep_Z)$. The group $G_S\cdot P_Z$ is compatible in $U^\ce$ and an orbit of $G_S\cdot P_Z$ is closed if and only if the corresponding orbit of $G$ is closed. Therefore we may assume $G=G_S\cdot P_Z$.

Dividing $P_Z$ by the ineffectivity of the $P_Z$-representation, we may furthermore assume that $P_Z$ acts effectively.
Since $\liep_Z$ acts by symmetric operators on $V$, we have a decomposition $V=\bigoplus_{i=1}^dV_i$ such that the action of the commutative group $P_Z$ on $V$ is given by a real character $\chi_i\colon P_Z\to\er^{>0}$ on each $V_i$. By assumption $P_Z$ acts effectively, so these characters define an injective homomorphism $\chi\colon P_Z\to(\er^{>0})^d$. Since $P_Z$ is contained in the center of $G$, each $V_i$ is invariant under the action of $G_S$. The universal complexification $G_S^\ce$ of $G_S$ in the sense of \cite{Hochschild} is complex reductive, contains $G_S$ as a compatible subgroup and by definition, it admits a holomorphic representation on $V_i^\ce$, which extends the $G_S$-representation. We define a representation of $(\ce^{*})^d$ on $V^\ce=\bigoplus_{i=1}^dV_i^\ce$ by letting the $i$-th component of $(\ce^*)^d$ act on $V_i^\ce$ by multiplication. The actions of $G_S^\ce$ and $(\ce^*)^d$ commute, so we obtain a holomorphic representation of the complex reductive group $G_S^\ce\times(\ce^*)^d$ on $V^\ce$. This group contains $G_S\times\chi(P_Z)$ as a compatible subgroup. Therefore we may assume $G=G_S\times\chi(P_Z)$,  $U^\ce=G_S^\ce\times(\ce^*)^d$ and $\hat V=V^\ce$. 

We define $H:=G_S\times (\er^{>0})^d$. Then $H$ is a compatible real form of $U^\ce$ which contains $G$. We denote its Cartan decomposition by $H=L\exp(\lieq)$. Note that by construction $L=K$ and that $\lieq$ is the direct sum of $\liep$ and a subspace $\liez$ of the center of $\lieu^\ce$. For an appropriate choice of the inner product, this sum can be assumed to be orthogonal. 

Altogether, for the proof of Theorem~\ref{theorem:closedOrbits}, we may assume that $G$ is a compatible subgroup of a compatible real form $H=K\exp(\lieq)$ of $U^\ce$ and that $H=G\times \exp(\liez)$.

With these assumptions, we can complete the proof of Theorem~\ref{theorem:closedOrbits}.

\begin{proof}[Proof of Theorem~\ref{theorem:closedOrbits}]
We prove the equivalent formulation of Theorem~\ref{theorem:closedOrbits} given in Section~\ref{section:rep}. Note that here $\pi(W)=\pw$ since $V$ is totally real in $V^\ce$. Again, let $n$ be maximal such that $\mathcal O_n\cap\pw\neq\emptyset$. Since $\mathcal O_n\cap\pw$ is open in $\pw$ with respect to the real Zariski topology and since $\mathcal S_G(\Mp)_c\cap\pw$ has non-empty interior in $\pw$ by assumption, there exists an $x\in\mathcal S_G(\Mp)_c\cap\mathcal O_n\cap\pw$. Let $g\in G$ with $gx\in\Mp$. Then $\beta:=\mu_\lieq(gx)$ is contained in the orthogonal complement of $\liep$ in $\lieq$ and in particular $\beta$ is contained in the center of $\lieu^\ce$.

By Proposition~\ref{proposition:Projektiv} the set $\mathcal S_H(\mathcal M_\lieq(\beta))_c\cap\mathcal O_n\cap\pw$ is Zariski open in $\pw$. It contains $x$ since $gx\in H\cdot x\cap\mathcal M_{\lieq}(\beta)$. Therefore it is non-empty. Note that it is $G$-invariant if $\pw$ is $G$-invariant. So for the proof of the theorem it now suffices to show that $\mathcal S_H(\Mq(\beta))_c\subset\mathcal S_G(\Mp)_c$.

For this, let $x\in \mathcal S_H(\Mq(\beta))_c$ and $y\in H\cdot x\cap\Mq(\beta)$. Then $y\in\Mp$ since $\beta$ is contained in the orthogonal complement of $\liep$. Consequently $G\cdot y$ is closed. The $G$-nullcone in $\hat V$ is $H$-invariant since $H=G\cdot \exp(\liez)$ where $\exp(\liez)$ is contained in the center of $U^\ce$. Therefore $\mathcal S_G(\Mp)$ is $H$-invariant by Lemma~\ref{lemma:semistabilZOffen} and we have $x\in H\cdot y\subset \mathcal S_G(\Mp)$. Let $h\in \exp(\liez)$ with $G\cdot x=G\cdot hy$. Then $G\cdot x=h\cdot G\cdot y$ is closed in $\mathcal S_G(\Mp)$ since $G\cdot y$ is closed.
\end{proof}

%%%%%%%%%%%%%%%%%%%%%%%%%%%%%%%%%%

\newcommand{\noopsort}[1]{} \newcommand{\printfirst}[2]{#1}
\newcommand{\singleletter}[1]{#1} \newcommand{\switchargs}[2]{#2#1}
\providecommand{\bysame}{\leavevmode\hbox to3em{\hrulefill}\thinspace}

\end{document}